\documentclass[a4paper,12pt]{article}
\usepackage{amssymb,amsmath,amsthm}
\usepackage{indentfirst}
\usepackage{amssymb,amsmath,amsthm}
\usepackage{euler}
\usepackage{hyperref,hypbmsec}
\newtheorem{Th}{\hskip\parindent Theorem}
\newtheorem{Le}{\hskip\parindent Lemma}

\newcommand{\N}{\mathbb{N}}

\newcounter{propet}

\renewcommand{\le}{\leqslant}\renewcommand{\ge}{\geqslant}

\begin{document}
\author{D.\,A.\,Frolenkov\footnote{The research was supported by the grant RFBR № 11-01-00759-a} }
\title{
Numerically explicit version of the P\'{o}lya--Vinogradov inequality
}
\date{}
\maketitle
\begin{abstract}
In this paper we proved a new numerically explicit version of the P\'{o}lya--Vinogradov inequality. Our proof is based on the new ideas of V.A. Bykovskii and improves a recent inequality obtained by C. Pomerance.
\end{abstract}
\vskip+1.4cm
In  paper ~\cite{Pomerance} C. Pomerance proved a new explicit version of the famous P\'{o}lya--Vinogradov inequality for character sums. In the present paper  we obtained an improvement of Pomerance's result. Our approach is based on a recent construction due to V. Bykovskii ~\cite{Bukovski}, ~\cite{Bukovski2}. In ~\cite{Bukovski}, ~\cite{Bukovski2} this construction was used to obtain new upper bounds for discrepancy of good lattice points sets. The present paper is organized as follows. In Section 1 we give a brief survey of classic and recent results on the topic. In Section 2 we formulate our main result. In Section 3  we formulate two lemmas by Pomerance. In Section 4 we describe Bykovskii's construction. In Section 5 we complete the proof of our main result.

\setcounter{Zam}0
\section{Introduction}
Let $\chi \pmod{q}$ be a primitive Dirichlet character. Put
$$S_{\chi}=\max_{0\le M<N\le q}\left|\sum_{n=M}^N\chi(n)\right|, \qquad T_{\chi}=\max_{N}\left|\sum_{a=0}^N\chi(a)\right|.$$
A character is defined to be even or odd if $\chi(-1)=1$ or $\chi(-1)=-1,$ respectively. In the case of even characters one has
\begin{gather}
S_{\chi}=2T_{\chi}.\label{T-S}
\end{gather}
In 1918 P\'{o}lya  ~\cite{Polya} and Vinogradov~\cite{Vinogradov} independently proved that for any nonprincipal Dirichlet character inequality
\begin{gather}
S_{\chi}\le c\sqrt{q}\log q\label{P-V}
\end{gather}
holds with an absolute constant c. In 2007 Granville and Soundararajan ~\cite{Granville-Sound} proved that for every primitive Dirichlet character $\chi \pmod{q}$ of odd order g
$$T_{\chi}\ll \sqrt{q}(\log q)^{1-\frac{\delta_g}{2}+o(1)}\quad \mbox{where}\quad \delta_g=1-\frac{g}{\pi}\sin\frac{\pi}{g},\;q\rightarrow\infty.$$
This result has been recently improved by Goldmakher see ~\cite{Gold}. He obtained the following result
$$T_{\chi}\ll \sqrt{q}(\log q)^{1-\delta_g+o(1)}.$$
Under the General Riemann Hypothesis Montgomery and Vaughan ~\cite[Theorem 3]{Mont-Vaughan} proved that
$$S_{\chi}\ll\sqrt{q}\log\log q .$$
In fact it is the best-possible result. Paley ~\cite{Paley} proved that there are an infinite class of quadratic characters $\chi_n \pmod{q_n}$ for which
$$S_{\chi_n}\gg \sqrt{q_n}\log\log q_n.$$
An important problem is to find the most precise form of inequality \eqref{P-V}.  There are two types of results.
The results of the first type do not take care about the explicit bounds for the remainder terms. The results of the second type give all constants explicitly. Usually the results of the first type have better constants in the main term.
\subsection{Asymptotic results}
Landau \cite{Landau} proved that
$$
S_{\chi}\le\left(\frac{1}{\pi\sqrt{2}}+o(1)\right)\sqrt{q}\log q \quad \mbox{if}\quad \chi(-1)=1
$$
and
$$
S_{\chi}\le\left(\frac{1}{2\pi}+o(1)\right)\sqrt{q}\log q \quad \mbox{if}\quad \chi(-1)=-1.
$$
Hildebrand \cite{Hildebrand1} obtained
$$
T_{\chi}\le\left(\frac{2}{3\pi^2}+o(1)\right)\sqrt{q}\log q \quad \mbox{if}\quad \chi(-1)=1
$$
and
$$
T_{\chi}\le\left(\frac{1}{3\pi}+o(1)\right)\sqrt{q}\log q \quad \mbox{if}\quad \chi(-1)=-1.
$$
Later Hildebrand \cite{Hildebrand2} improved his result for even characters. He showed that the estimate
$$
T_{\chi}\le\left(\frac{c}{\pi\sqrt{3}}+o(1)\right)\sqrt{q}\log q \quad \mbox{if}\quad \chi(-1)=1
$$
where
$$
c=
\left\{
              \begin{array}{ll}
                \frac{1}{4}, & \hbox{if q is cubefree;} \\
                \frac{1}{3}, & \hbox{else,}
               \end{array}
\right.
$$
holds.
Granville and Soundararajan \cite{Granville-Sound} obtained two inequalities
$$
T_{\chi}\le\left(\frac{69}{70}\frac{c}{\pi\sqrt{3}}+o(1)\right)\sqrt{q}\log q \quad \mbox{if}\quad \chi(-1)=1
$$
and
$$
T_{\chi}\le\left(\frac{c}{\pi}+o(1)\right)\sqrt{q}\log q \quad \mbox{if}\quad \chi(-1)=-1.
$$
Up to now this result is the best-known one.
\subsection{Numerically explicit results}
In this section we discus numerically explicit versions of the P\'{o}lya--Vinogradov inequality.
Qiu \cite{Qiu} proved that
$$
S_{\chi}\le\frac{4}{\pi^2}\sqrt{q}\log q+0.38\sqrt{q}+0.608\frac{1}{\sqrt{q}}+0.116\sqrt{q}.
$$
Simalarides \cite{Simalarides}, \cite{Simalarides2} obtained estimates
$$
T_{\chi}\le\frac{3}{4\pi}\sqrt{q}\log q+\left(2-\frac{\log 2}{\pi}-\frac{\gamma}{2\pi}\right) \quad \mbox{if}\quad \chi(-1)=1
$$
and
$$
T_{\chi}\le\frac{1}{\pi}\sqrt{q}\log q+\sqrt{q}+\frac{1}{2} \quad \mbox{if}\quad \chi(-1)=-1.
$$
Dobrowolski  and Williams \cite{Dobro-Williams} proved that for any nonprincipal Dirichlet character $\chi \pmod{q}$ one has
$$
S_{\chi}\le\frac{1}{2\log2}\sqrt{q}\log q+3\sqrt{q}.
$$
Bachman and Rachakonda \cite{Bachman-Racha} improved their result and proved
$$
S_{\chi}\le\frac{1}{3\log3}\sqrt{q}\log q+6.5\sqrt{q}
$$
for any nonprincipal Dirichlet character $\chi \pmod{q}.$
In ~\cite{Pomerance}  Pomerance proved that
\begin{gather}\label{P1}
S_{\chi}\le\frac{2}{\pi^2}\sqrt{q}\log q+\frac{4}{\pi^2}\sqrt{q}\log\log q+\frac{3}{2}\sqrt{q}\quad \mbox{if}\quad \chi(-1)=1
\end{gather}
and
\begin{gather}\label{P2}
S_{\chi}\le\frac{1}{2\pi}\sqrt{q}\log q+\frac{1}{\pi}\sqrt{q}\log\log q+\sqrt{q}\quad \mbox{if}\quad \chi(-1)=-1.
\end{gather}
Up to now these bounds are the best-known numerically explicit versions of the P\'{o}lya--Vinogradov inequality.
\section{Main result}
We prove the following theorem which improves both \eqref{P1} and \eqref{P2}  in the second term.
\begin{Th}\label{Th1}
Let $\chi \pmod{q}$ be a primitive character then
\begin{enumerate}
  \item If $\chi(-1)=1$ then
$$
S_{\chi}\le\frac{2}{\pi^2}\sqrt{q}\log q+\frac{4}{\pi^2}\sqrt{q}
\left(1+\gamma+\log C_0\right)+\psi_1(q),
$$
  \item If $\chi(-1)=-1$ then
$$
S_{\chi}\le\frac{1}{2\pi}\sqrt{q}\log q+\frac{1}{\pi}\sqrt{q}\left(1+\gamma+\log\frac{2C_0}{\pi}\right)+\psi_2(q),
$$
where $\gamma$ is the Euler constant, $C_0=4\pi^{5/2}+5$ and
$$\psi_1(q)=1+\frac{24}{\pi^2C_0}+\frac{8}{\pi^2}\frac{\sqrt{q}}{exp(\frac{2\sqrt{q}}{C_0})-1}$$
$$\psi_2(q)=1+\frac{3}{C_0}+\frac{2}{\pi}\frac{\sqrt{q}}{exp(\frac{\pi\sqrt{q}}{C_0})-1}$$
\end{enumerate}
\end{Th}
Easy calculations show that for $q>18000$ the first case of Theorem 1 is better than \eqref{P1} and
for $q>28000$ the second case of Theorem 1 is better than \eqref{P2}.
The proof of this theorem is based on the ideas of Bykovskii ~\cite{Bukovski} and on the ideas of Pomerance~\cite{Pomerance}.
\section{Pomerance's lemmas}
In this section we formulate all necessary results from the Pomerance's paper~\cite{Pomerance}.
\begin{Le}\label{sin}
For all real numbers $x$ and positive integers $n$, we have
$$
\sum_{j=1}^{n}\frac{|\sin jx|}{j}<\frac{2}{\pi}\log n+\frac{2}{\pi}\left(\gamma+\log2+\frac{3}{n}\right).
$$
\begin{proof}
See ~\cite[Lemma 3]{Pomerance}.
\end{proof}
\end{Le}
\begin{Le}\label{cos}
For all real numbers $\alpha$, $\beta$ and positive integers $n$, we have
$$
\sum_{m=1}^{n}\frac{1}{m}\left|\cos m\alpha-\cos m\beta\right|<\log n+\gamma+\log2+\frac{3}{n}.
$$
\begin{proof}
See ~\cite[\S 4]{Pomerance}.
\end{proof}
\end{Le}

\section{Bykovskii's method}
In this section we describe the main construction from the paper ~\cite{Bukovski}. We give a slight modification. Let $\theta:[0,\infty)\rightarrow[0,1]$ be the following function
\begin{gather*}
\theta(t)=
\left\{
              \begin{array}{ll}
                0, & \hbox{if $0<t\le\frac{1}{2}$;} \\
                2t-1, & \hbox{if $\frac{1}{2}\le t\le1$;} \\
                2-t, & \hbox{if $1\le t\le2$;} \\
                0, & \hbox{if $2\le t$.}
              \end{array}
\right.
\end{gather*}
For all $t\in(0,\infty)$ one has
$$\sum_{j=-\infty}^{\infty}\theta\left(\frac{t}{2^j}\right)=1.$$
Let $\omega,\omega^{*}:(0,\infty)\rightarrow[0,1]$be the following functions
\begin{gather*}
\omega(t)=\sum_{j=0}^{\infty}\theta\left(\frac{t}{2^j}\right)=
\left\{
              \begin{array}{ll}
                0, & \hbox{if $0<t\le\frac{1}{2}$;} \\
                2t-1, & \hbox{if $\frac{1}{2}\le t\le1$;} \\
                1, & \hbox{if $1\le t$,}
              \end{array}
\right.
\end{gather*}
\begin{gather*}
\omega^{*}(t)=1-\omega(t)=\sum_{j=-\infty}^{-1}\theta\left(\frac{t}{2^j}\right)=
\left\{
              \begin{array}{ll}
                1, & \hbox{if $0<t\le\frac{1}{2}$;} \\
                2-2t, & \hbox{if $\frac{1}{2}\le t\le1$;} \\
                0, & \hbox{if $1\le t$.}
              \end{array}
\right.
\end{gather*}
For any positive numbers $P$ and $P'$ put
$$
S(u;P,P')=\sum_{m=1}^{\infty}\frac{1}{m}\omega\left(\frac{m}{P}\right)\omega^{*}\left(\frac{m}{P'}\right)\sin 2\pi mu.
$$
If $2P\le P'$ then
$$S(u;P,P')=S(u;P)-S(u;P'),$$
where
\begin{gather}\label{SuP}
S(u;P)=\sum_{m=1}^{\infty}\frac{1}{m}\omega\left(\frac{m}{P}\right)\sin 2\pi mu.
\end{gather}
For any positive number $P$ let
\begin{gather}\label{Gpu}
G_{P}(u)=\sum_{n=-\infty}^{\infty}\frac{1}{1+P^2(n+u)^2}=
\frac{\pi}{P}\sum_{n=-\infty}^{\infty}exp\left(-\frac{|n|}{P}\right)e(nu),
\end{gather}
where $e(u)=\exp(2\pi iu).$
Bykovskii proved (see~\cite[Lemma 1]{Bukovski}) that
$$\left|S(u;P)\right|\le 14G_{P}(u).$$
We improve this result.
\begin{Le}\label{SuP<Gpu}
For any $u\in [0;1]$ and $P>1$ the inequality
$$\left|S(u;P)\right|\le \frac{C_0}{2\pi^2}G_{P}(u)$$
is valid.
\begin{proof}
Take $A\in [0;1]$. The optimal value of parameter $A$ will be calculated later. We should consider two cases.
\begin{enumerate}
  \item Let $P\|u\|\le A.$ \par
  We see that
\begin{gather}\label{sin}
|\sin 2\pi mu|\le 2\pi\|mu\|
\end{gather}
and
\begin{gather}\label{sumsin}
\left|\sum_{m=m_1}^{m_2}\sin 2\pi mu\right|\le\frac{1}{2\|u\|}.
\end{gather}
The proof of \eqref{sumsin} can be found in ~\cite[Lemma 1]{Bukovski}. By partial summation and \eqref{sumsin} we have
\begin{gather}\label{sin/n}
\left|\sum_{m>\frac{B}{\|u\|}}\frac{\sin 2\pi mu}{m}\right|\le \frac{1}{2\|u\|}\frac{\|u\|}{B}\le\frac{1}{2B}.
\end{gather}
Let $B\in [A;1]$. The optimal value of parameter $B$ will be calculated later.We will define $B$ later. We see that
$$P\le\frac{A}{\|u\|}\le\frac{B}{\|u\|}.$$
So by \eqref{SuP},\;\eqref{sumsin},\;\eqref{sin/n} we have
\begin{gather*}
\left|S(u;P)\right|\le\sum_{m\le\frac{B}{\|u\|}}\frac{|\sin 2\pi mu|}{m}+\left|\sum_{m>\frac{B}{\|u\|}}\frac{\sin 2\pi mu}{m}\right|\le \sum_{m\le\frac{B}{\|u\|}}\frac{2\pi\|mu\|}{m}+\frac{1}{2B}.
\end{gather*}
If $B\le\frac{1}{2}$ then
$$\left|S(u;P)\right|\le2\pi B+\frac{1}{2B},$$
else
\begin{gather*}
\left|S(u;P)\right|\le 2\pi\frac{1}{2}+2\pi\sum_{\frac{1}{2\|u\|}<m\le\frac{B}{\|u\|}}\frac{1-mu}{m}+\frac{1}{2B}\le
2\pi\left(\frac{1}{2}+\log 2B-(B-\frac{1}{2})\right)+\frac{1}{2B}=\\=2\pi\left(1-B+\log 2B\right)+\frac{1}{2B}.
\end{gather*}
By the definition of function $G_{P}(u)$ we have
$$G_{P}(u)\ge\frac{1}{1+P^2\|u\|^2}\ge\frac{1}{1+A^2}.$$
So
$$\left|S(u;P)\right|\le C_1(1+A^2)G_{P}(u),$$
where
\begin{gather*}
C_1=
\left\{
              \begin{array}{ll}
                2\pi B+\frac{1}{2B}, & \hbox{if $A\le B\le\frac{1}{2}$;} \\
                2\pi\left(1-B+\log 2B\right)+\frac{1}{2B}, & \hbox{if $\max\left(\frac{1}{2},A\right)\le B\le1$.}
               \end{array}
\right.
\end{gather*}
  \item Let $P\|u\|> A.$\par
Our proof has much in common with the proof of the second case in~\cite[Lemma 1]{Bukovski}. So we use the notation of ~\cite[Lemma 1]{Bukovski}.  It was shown in ~\cite[Lemma 1]{Bukovski} that for any $P'>2P$ we have
\begin{gather*}
\left|S(u;P,P')\right|=\frac{1}{2}\left|\sum_{n=-\infty}^{\infty}(\eta(u-\|u\|;P,P')-\eta(u+\|u\|;P,P'))\right|,
\end{gather*}
where
$$\eta(w;P,P')=\int_{P/2}^{P}\omega\left(\frac{x}{P}\right)\omega^{*}\left(\frac{x}{P'}\right)e(-wx)\frac{dx}{x}.$$
By the ideas of ~\cite[Lemma 1]{Bukovski} we have
$$
|\eta(w;P,P')|\le\frac{1}{4\pi^2w^2}\left(\frac{10}{P^2}+\frac{16}{P'^2}\right).
$$
So
$$
\left|S(u;P,P')\right|\le\sum_{n=-\infty}^{\infty}\frac{1}{4\pi^2(n+u)^2}\left(\frac{10}{P^2}+\frac{16}{P'^2}\right).
$$
Hence
$$
\left|S(u;P)\right|\le\frac{5}{2\pi^2}\sum_{n=-\infty}^{\infty}\frac{1}{P^2(n+u)^2}.
$$
For any integer number $n$ we have
$$
P^2(n+u)^2>\frac{A^2}{1+A^2}\left(1+P^2(n+u)^2\right).
$$
Therefore
$$
\left|S(u;P)\right|\le\frac{5}{2\pi^2}\left(1+\frac{1}{A^2}\right)\sum_{n=-\infty}^{\infty}\frac{1}{1+P^2(n+u)^2}<
\frac{5}{2\pi^2}\left(1+\frac{1}{A^2}\right)G_P(u).
$$
\end{enumerate}
It is easy to prove that
$$f_1(A)=C_1(1+A^2)$$
is an increasing function when $A\in [0;1]$. So we should take $A$ such that
$$f_1(A)=\frac{5}{2\pi^2}\left(1+\frac{1}{A^2}\right).$$
If $A\le\frac{1}{2\sqrt{\pi}}$ then
$$
\min_{A\le B\le\frac{1}{2}}\left(2\pi B+\frac{1}{2B}\right)=2\sqrt{\pi}.
$$
Hence
$$
f_1(A)=2\sqrt{\pi}(1+A^2), \quad 0\le A\le\frac{1}{2\sqrt{\pi}}.
$$
Solving the system
\begin{gather}
\left\{
  \begin{array}{ll}
    2\sqrt{\pi}(1+A^2)=\frac{5}{2\pi^2}\left(1+\frac{1}{A^2}\right), &  \\
    0\le A\le\frac{1}{2\sqrt{\pi}},
   \end{array}
\right.
\end{gather}
we have $$A=\frac{\sqrt{5}}{2\pi^{5/4}}.$$ So
$$
\left|S(u;P)\right|\le f_1\left(\frac{\sqrt{5}}{2\pi^{5/4}}\right)G_{P}(u)=\left(2\sqrt{\pi}+\frac{5}{2\pi^2}\right)G_{P}(u).
$$
Lemma is proved.
\end{proof}
\end{Le}

\begin{Le}\label{Gpu-sum}
For any $q\in\N$ and $P>1$ the formula
$$
\sum_{a=1}^{q}G_{P}\left(\frac{a}{q}\right)=\pi\frac{q}{P}+2\pi\frac{q}{P}\frac{1}{\exp{\left(\frac{q}{P}\right)-1}}
$$
is valid.
\begin{proof}
Set
\begin{gather}\label{deltasum}
\delta_q(a)=\frac{1}{q}\sum_{x=1}^q e\left(\frac{ax}{q}\right)=
\left\{
              \begin{array}{ll}
                1, & \hbox{if $a\equiv 0 \pmod{q}$;} \\
                0, & \hbox{else,}
              \end{array}
\right.
\end{gather}
then by \eqref{Gpu} we have
\begin{gather*}
\sum_{a=1}^{q}G_{P}\left(\frac{a}{q}\right)=\frac{\pi}{P}\sum_{a=1}^{q}\sum_{n=-\infty}^{\infty}
\exp\left(-\frac{|n|}{P}\right)e\left(n\frac{a}{q}\right)=
\frac{\pi}{P}\sum_{n=-\infty}^{\infty}\exp\left(-\frac{|n|}{P}\right)q\delta_q(n)=\\=
\pi\frac{q}{P}\left(1+2\sum_{n=1}^{\infty}\exp\left(-\frac{n}{P}\right)\delta_q(n)\right)=
\pi\frac{q}{P}\left(1+2\sum_{n=1}^{\infty}\exp\left(-\frac{qn}{P}\right)\right)=
\pi\frac{q}{P}\left(1+2\frac{\exp\left(-\frac{q}{P}\right)}{1-\exp\left(-\frac{q}{P}\right)}\right).
\end{gather*}
Lemma is proved.
\end{proof}
\end{Le}

\section{Proof of Theorem \ref{Th1}}
Let $0\le a<b\le1$ and $\lambda(x;a,b)$ be the function on $[0;1)$ satisfying
\begin{gather*}
\lambda(x;a,b)=
\left\{
              \begin{array}{ll}
                \frac{1}{2}, & \hbox{if $x=a$;} \\
                1, & \hbox{if $a<x<b$;} \\
                \frac{1}{2}, & \hbox{if $x=b$;}\\
                0, & \hbox{else.}
               \end{array}
\right.
\end{gather*}
Writing its Fourier expansion we have
$$
\lambda(x;a,b)=b-a+\frac{1}{\pi}S(x-a)-\frac{1}{\pi}S(x-b)
$$
where
$$
S(x)=\sum_{m=1}^{\infty}\frac{\sin2\pi mx}{m}.
$$
For any $0\le M<N\le q$ one has
\begin{gather}
\sum_{a=M}^N\chi(a)=\sum_{a=1}^q\chi(a)\lambda\left(a;\frac{M}{q},\frac{N}{q}\right)+\frac{\chi(M)+\chi(N)}{2}=\notag\\=
\frac{1}{\pi}\sum_{a=1}^q\chi(a)\left(S\left(\frac{a-M}{q}\right)-S\left(\frac{a-N}{q}\right)\right)+\frac{\chi(M)+\chi(N)}{2}.\label{fourier}
\end{gather}
The Gauss sum $\tau(\chi)$ is defined as
$$
\tau(\chi)=\sum_{a=1}^q\chi(a)e\left(\frac{a}{q}\right).
$$
For a primitive Dirichlet character $\chi \pmod{q}$ one has
\begin{gather}
|\tau(\chi)|=\sqrt{q}.\label{gaussmodul}
\end{gather}
It was proved in ~\cite[(10)]{Pomerance} that
\begin{gather}
\overline{\chi}(m)\tau(\chi)=
\left\{
              \begin{array}{ll}
              \sum\limits_{a=1}^q\chi(a)\cos\frac{2\pi am}{q}, & \hbox{if $\chi$ is even;} \\
              i\sum\limits_{a=1}^q\chi(a)\sin\frac{2\pi am}{q}, & \hbox{if $\chi$ is odd.} \\
               \end{array}
\right.\label{gauss}
\end{gather}

\subsection{The case of even characters}
It follows from \eqref{fourier} that
\begin{gather*}
\sum_{a=1}^N\chi(a)=
\frac{1}{\pi}\sum_{a=1}^q\chi(a)\left(S\left(\frac{a}{q}\right)-S\left(\frac{a-N}{q}\right)\right)+\frac{\chi(N)}{2}.
\end{gather*}
Let $P>0$ by the definition of functions  $\omega(t),\omega^{*}(t)$ we have
\begin{gather}
S(u)=\sum_{m\le P}\frac{\sin 2\pi mu}{m}\omega^{*}\left(\frac{m}{P}\right)+\sum_{m=1}^{\infty}\frac{\sin 2\pi mu}{m}\omega\left(\frac{m}{P}\right)=\notag\\=
\sum_{m\le P}\frac{\sin 2\pi mu}{m}\omega^{*}\left(\frac{m}{P}\right)+S(u;P).\label{S-SP}
\end{gather}
Hence
\begin{gather}
\sum_{a=1}^N\chi(a)=
\frac{1}{\pi}\sum_{a=1}^q\chi(a)\sum_{m\le P}\frac{\sin 2\pi m\frac{a}{q}-
\sin 2\pi m\frac{a-N}{q}}{m}\omega^{*}\left(\frac{m}{P}\right)+\notag\\+
\frac{1}{\pi}\sum_{a=1}^q\chi(a)\left(S\left(\frac{a}{q};P\right)-S\left(\frac{a-N}{q};P\right)
\right)
+\frac{\chi(N)}{2}.\label{evenest}
\end{gather}
By \eqref{gauss} one has
\begin{gather*}
\frac{1}{\pi}\sum_{a=1}^q\chi(a)\sum_{m\le P}\frac{\sin 2\pi m\frac{a}{q}-
\sin 2\pi m\frac{a-N}{q}}{m}\omega^{*}\left(\frac{m}{P}\right)=
\frac{1}{\pi}\sum_{m\le P}\frac{\overline{\chi}(m)\tau(\chi)}{m}\sin 2\pi m\frac{N}{q}\omega^{*}\left(\frac{m}{P}\right).
\end{gather*}
By  \eqref{gaussmodul} one has
\begin{gather}
\left|\frac{1}{\pi}\sum_{a=1}^q\chi(a)\sum_{m\le P}\frac{\sin 2\pi m\frac{a}{q}-
\sin 2\pi m\frac{a-N}{q}}{m}\omega^{*}\left(\frac{m}{P}\right)\right|\le\frac{\sqrt{q}}{\pi}\sum_{m\le P}
\frac{\left|\sin 2\pi m\frac{N}{q}\right|}{m}\label{evensin}.
\end{gather}
Thus, by \eqref{evensin} and by Lemma \ref{SuP<Gpu},
\begin{gather}
\left|\sum_{a=1}^N\chi(a)\right|\le\frac{\sqrt{q}}{\pi}\sum_{m\le P}
\frac{\left|\sin 2\pi m\frac{N}{q}\right|}{m}+\frac{2}{\pi}\sum_{a=1}^q\frac{C_0}{2\pi^2}G_{P}\left(\frac{a}{q}\right)
+\frac{1}{2}.
\end{gather}
By Lemma \ref{sin} and Lemma \ref{Gpu-sum} we have
\begin{gather}
\left|\sum_{a=1}^N\chi(a)\right|\le
\frac{2}{\pi^2}\sqrt{q}\left(\log P+\gamma+\log2+\frac{3}{P}\right)+
2\frac{C_0}{2\pi^2}\frac{q}{P}+2\frac{C_0}{2\pi^2}\frac{q}{P}\frac{2}{\exp{\left(\frac{q}{P}\right)-1}}
+\frac{1}{2}.
\end{gather}
Now we take $P=\frac{C_0}{2}\sqrt{q}$. So
\begin{gather}
\left|\sum_{a=1}^N\chi(a)\right|\le
\frac{1}{\pi^2}\sqrt{q}\log q+\frac{2}{\pi^2}\sqrt{q}
\left(1+\gamma+\log C_0\right)+\frac{12}{\pi^2C_0}+\frac{4}{\pi^2}\frac{\sqrt{q}}{exp(\frac{2\sqrt{q}}{C_0})-1}
+\frac{1}{2}.
\end{gather}
We take into account \eqref{T-S} to complete the proof of Theorem \ref{Th1} in the case of even characters.
\subsection{The case of odd characters}
By \eqref{fourier} and \eqref{S-SP} we have
\begin{gather}
\sum_{a=M}^N\chi(a)=
\frac{1}{\pi}\sum_{a=1}^q\chi(a)\sum_{m\le P}\frac{\sin 2\pi m\frac{a-M}{q}-
\sin 2\pi m\frac{a-N}{q}}{m}\omega^{*}\left(\frac{m}{P}\right)+\notag\\+
\frac{1}{\pi}\sum_{a=1}^q\chi(a)\left(S\left(\frac{a-M}{q};P\right)-S\left(\frac{a-N}{q};P\right)
\right)
+\frac{\chi(M)+\chi(N)}{2}.\label{oddsum}
\end{gather}
Using  \eqref{gauss}
\begin{gather*}
\frac{1}{\pi}\sum_{a=1}^q\chi(a)\sum_{m\le P}\frac{\sin 2\pi m\frac{a-M}{q}-
\sin 2\pi m\frac{a-N}{q}}{m}\omega^{*}\left(\frac{m}{P}\right)=\\=-
\frac{i}{\pi}\sum_{m\le P}\frac{\overline{\chi}(m)\tau(\chi)}{m}
\left(\cos2\pi m\frac{M}{q}-\cos 2\pi m\frac{N}{q}\right)\omega^{*}\left(\frac{m}{P}\right)
\end{gather*}
Using  \eqref{gaussmodul} and Lemma \ref{cos} we have
\begin{gather}
\left|\frac{1}{\pi}\sum_{a=1}^q\chi(a)\sum_{m\le P}\frac{\sin 2\pi m\frac{a-M}{q}-
\sin 2\pi m\frac{a-N}{q}}{m}\omega^{*}\left(\frac{m}{P}\right)\right|\le\frac{\sqrt{q}}{\pi}\sum_{m\le P}
\frac{1}{m}\left|\cos2\pi m\frac{M}{q}-\cos 2\pi m\frac{N}{q}\right|\le\notag\\\le
\frac{\sqrt{q}}{\pi}\left(\log P+\gamma+\log2+\frac{3}{P}\right).\label{oddcos}
\end{gather}
Thus, by \eqref{oddsum}, \eqref{oddcos} and by Lemma \ref{SuP<Gpu} we have
\begin{gather*}
\left|\sum_{a=M}^N\chi(a)\right|\le\frac{\sqrt{q}}{\pi}\left(\log P+\gamma+\log2+\frac{3}{P}\right)+
\frac{2}{\pi}\sum_{a=1}^q\frac{C_0}{2\pi^2}G_{P}\left(\frac{a}{q}\right)+1
\end{gather*}
By Lemma \ref{Gpu-sum} we have
\begin{gather*}
\left|\sum_{a=M}^N\chi(a)\right|\le\frac{\sqrt{q}}{\pi}\left(\log P+\gamma+\log2+\frac{3}{P}\right)+
2\frac{C_0}{2\pi^2}\frac{q}{P}+2\frac{C_0}{2\pi^2}\frac{q}{P}\frac{2}{\exp{\left(\frac{q}{P}\right)-1}}+1
\end{gather*}
Now we take $P=\frac{C_0}{\pi}\sqrt{q}$. So
\begin{gather*}
\left|\sum_{a=M}^N\chi(a)\right|\le\frac{1}{2\pi}\sqrt{q}\log q+\frac{1}{\pi}\sqrt{q}\left(1+\gamma+\log\frac{2C_0}{\pi}\right)+\frac{3}{C_0}+
\frac{2}{\pi}\frac{\sqrt{q}}{exp(\frac{\pi\sqrt{q}}{C_0})-1}
+1.
\end{gather*}
Theorem \ref{Th1} is proved.

\textit{D.A. Frolenkov}\\
\textit{Department of Number theory}\\
\textit{Moscow State University}\\
\textit{e-mail: frolenkov\underline{ }adv@mail.ru}


\newpage
\begin{thebibliography}{99}
\bibitem{Polya} \textsc{P\'{o}lya ~G.}
{\"{U}ber die Verteilung der quadratischen Reste und Nichreste. Nachrichten K\"{o}nigl.Ges.Wiss.G\"{o}ttingen (1918),pp. 21-29.}
\bibitem{Vinogradov} \textsc{Vinogradov ~I.\,M.}
{On the distribution of power residues and non-residues, J. Phys. Math. Soc. Perm. Univ. 1 (1918),pp. 94-98; Selected works, Springer Berlin,1985,pp. 53-56.}
\bibitem{Landau} \textsc{Landau E.}
{Absch\"{a}tzungen von Charaktersummen, Einheiten und klassenzahlen. Nachrichten K\"{o}nigl.Ges.Wiss.G\"{o}ttingen (1918),pp. 79-97.}
\bibitem{Hildebrand1} \textsc{Hildebrand ~A.}
{On the constant in the P\'{o}lya--Vinogradov inequality, Canad. Math. Bull. 31 (1988), 347-352.}
\bibitem{Hildebrand2} \textsc{Hildebrand ~A.}
{Large values of character sums, J. Number Theory 29 (1988), 271-296.}
\bibitem{Qiu} \textsc{Qiu ~Z.\,M.}
{An inequality of Vinogradov for character (Chinese), J. Shandong Univ., Nat. Sci. Ed. 26 (1991), 125–128.}
\bibitem{Simalarides} \textsc{Simalarides ~A.\,D.}
{An elementary proof of P\'{o}lya--Vinogradov's inequality. Period. Math.Hungar. 38 (1999) 99-101.}
\bibitem{Simalarides2} \textsc{Simalarides ~A.\,D.}
{An elementary proof of P\'{o}lya--Vinogradov's inequality, 2. Period. Math.Hungar. 40 (2000) 71-75.}
\bibitem{Dobro-Williams} \textsc{Dobrowolski ~E. and Williams ~K.\,S.}
{An upper bound for the sum $\sum_{n=a+1}^{a+H}f(n)$ for a certain class of functions f, Proc. Amer. Math. Soc. 114(1992), 29-35}
\bibitem{Bachman-Racha} \textsc{Bachman ~G. and Rachakonda ~L.}
{On a problem of Dobrowolski and Williams and the P\'{o}lya--Vinogradov inequality, Ramanujan J. 5 (2001), 65-71.}
\bibitem{Pomerance} \textsc{Pomerance ~C.}
{Remarks on the P\'{o}lya--Vinogradov inequality. Integers (Proceedings of the Integers Conference, October 2009), 11A (2011), Article 19, 11pp.}
\bibitem{Bukovski} \textsc{Bykovskii ~V.\,A.}
{The disperancy of the Korobov lattice points. To appear in Izv. RAN. Ser. Mat.}
\bibitem{Bukovski2} \textsc{Bykovskii ~V.\,A.}
{The discrepancy of the Korobov lattice points. Program and Abstrct Book p.10-11. 27th Journ\'{e}es Arithm\'{e}tiques confernce.}
\bibitem{Paley} \textsc{Paley ~R.\,E.\,A.\,C.}
{A theorem on characters, J.London Math. Soc. 7 (1932), 28-32. }
\bibitem{Mont-Vaughan} \textsc{Monjtgomery ~H.\,L. and Vaughan ~R.\,C.}
{Exponential sums with multiplicative coefficients, Invent. Math. 43 (1977), 69-82. }
\bibitem{Granville-Sound} \textsc{Granville ~A. and Soundararajan ~K.}
{Large character sums: pretentious characters and the  P\'{o}lya--Vinogradov theorem, Jour. AMS Vol. 20, Number 2(2007), 357-384. }
\bibitem{Gold} \textsc{Goldmakher ~L.}
{Multiplicative mimicry and improvements of the P\'{o}lya-Vinogradov inequality by Goldmakher, Leo I., Ph.D., UNIVERSITY OF MICHIGAN, 2009, 109 pages; 3382190, preprint is available in  arXiv:0911.5547v2}
\end{thebibliography}
\end{document}